\documentclass[12pt]{amsart}
\usepackage{graphicx,epsf,afterpage,url}
\usepackage{amsfonts,amsmath,amssymb,amscd,amsthm}
\input{xy}
\xyoption{all}
\def\Z{{\mathbb Z}} \def\R{{\mathbb R}}

\def\t{\widetilde}

\long\def\comment#1\endcomment{}

\graphicspath{{morefig/}{figures/}{figures00/}{figures02/}{figures05/}{figlms/}{obstruct-final/}{0808.1187v1/}}

\usepackage[usenames,dvipsnames]{color}

\theoremstyle{theorem}
\newtheorem{Theorem}{Theorem}[section]
\newtheorem{theorem}[Theorem]{Theorem}

\newtheorem{corollary}[Theorem]{Corollary}
\newtheorem{proposition}[Theorem]{Proposition}
\newtheorem{conjecture}[Theorem]{Conjecture}
\newtheorem{remark}[Theorem]{Remark}
\newtheorem{example}[Theorem]{Example}
\newtheorem{problem}[Theorem]{Problem}

\theoremstyle{definition}

\begin{document}

%\begin{frame}

\title{Stability of intersections of graphs in the plane\\
and the van Kampen obstruction}

%\smallskip\centerline{\bf
%{\bf Acknowledgements.}

\author{Arkadiy Skopenkov}

\address{Moscow Institute of Physics and Technology and Independent University of Moscow.
Info: \url{www.mccme.ru/~skopenko}}

\email{skopenko@mccme.ru}


%derivative of a path

\date{}

\maketitle

\begin{abstract}
A map $\varphi:K\to \R^2$ of a graph $K$ is {\it approximable by embeddings}, if for each $\varepsilon>0$ there is an $\varepsilon$-close to $\varphi$ embedding $f:K\to \R^2$.
Analogous notions were studied in computer science under the names of {\it cluster planarity} and {\it weak simplicity}.
This short survey is intended not only for specialists in the area, but also for mathematicians from other areas.

We present criteria for approximability by embeddings (P. Minc, 1997, M. Skopenkov, 2003) and their algorithmic corollaries.
We introduce {\it the van Kampen (or Hanani-Tutte) obstruction} for approximability by embeddings and discuss its completeness. We discuss analogous problems of moving graphs in the plane apart (cf. S.~Spie\. z and H.~Toru\'nczyk, 1991) and finding closest embeddings (H. Edelsbrunner).
We present higher dimensional generalizations, including completeness of the van Kampen obstruction and its algorithmic corollary (D. Repov\v s and A. Skopenkov, 1998).
\end{abstract}

\tableofcontents

%A map $\varphi:K\sqcup L\to \R^2$ of the disjoint union of graphs $K$ and $L$ is {\bf approximable by maps
%with disjoint images}, if for each $\varepsilon>0$ there is an $\varepsilon$-close to $\varphi$ map
%$f:K\sqcup L\to\R^2$ such that $f(K)\cap f(L)=\emptyset$. We present open problems on this notion.
% $r$-tuple intersections.
%and of the triviality of the $r$-tuple obstruction
%for $r$ not a prime power
%(M. \"Ozaydin, 1987, unpublished).

%\end{frame}

%\begin{frame}

\section{Approximability by embeddings}

\subsection{Definition, examples and discussion}

All the maps below are tacitly assumed to be continuous or piecewise-linear (PL).

A map $\varphi:K\to \R^2$ of a graph $K$ is {\bf approximable by embeddings}, if for each $\varepsilon>0$ there is an $\varepsilon$-close to $\varphi$ embedding $f:K\to\R^2$.

Even the cases when $\varphi$ is either a path or a cycle, i. e. either $K\cong I:=[0,1]$ or
$K\cong S^1:=\{(x,y)\in\R^2\ :\ x^2+y^2=1\}$, is interesting.

This notion appeared in studies of planarity of compacta and dynamical systems, see
%later
some illustrating examples in \S\ref{s:repl}, \S\ref{s:exdy}.
Related notions of {\it cluster planarity, weak simplicity, projected embedding, level-planarity and
monotone drawings} appeared in computer science and differential topology,
see \cite{ARS, FPSS, Fu14, PT}, \cite[\S3.5]{Sk}, \cite[Theorem 1]{Fu16} and references therein.

{\it The following examples} (see also \cite[Fig. 5]{RS}) show that this notion is interesting in itself.
See elementary introduction in \cite{DNSS, L}.

%\pause

\begin{figure}[h]\centering
\includegraphics{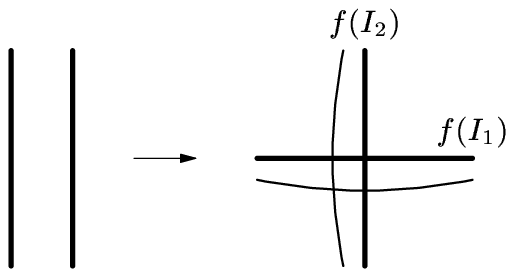}\qquad \qquad
\includegraphics{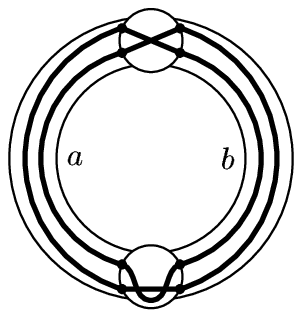}
\caption{Transversal intersection and the standard 2-winding are not approximable by embeddings}\label{hund}
\end{figure}

A {\it transversal self-intersection} of a PL map $\varphi:K\to \R^2$ is a pair of disjoint arcs $i,j\subset K$
such that $\varphi i$ and $\varphi j$ intersect transversally in the plane.

\begin{figure}[h]\centering
\begin{tabular}{@{}ccc@{}}
\includegraphics{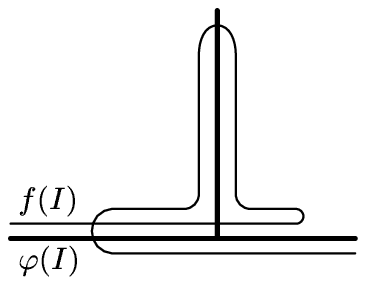} \qquad & \includegraphics{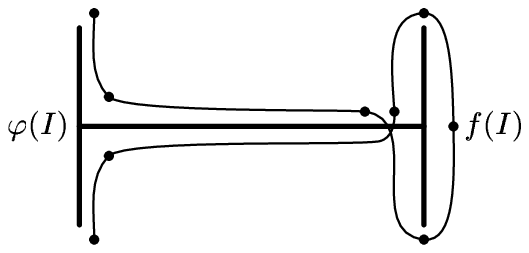} & \qquad \includegraphics{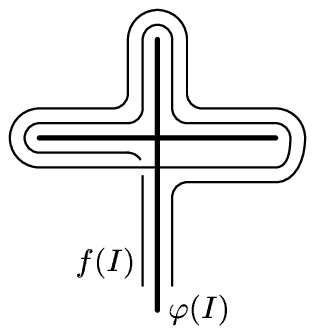}\\
(a) & (b) & (c)
\end{tabular}
\caption{Paths not approximable by embeddings and having no transversal intersections}\label{putibcd}
\end{figure}

\begin{remark}
(a) A map from a graph to a point in the plane is approximable by embeddings if and only if the graph is planar.

(b) Any map $f:I\to I\subset\R^2$ is approximable by embeddings.
(Indeed, take the graph of $f$ in $I\times I\subset\R^2$ and compress it to the first factor.)
\end{remark}

\begin{proposition}
(a) Any map $T\to I\subset\R^2$ is approximable by embeddings, where $T$ is triod. (Sieklucki, 1969 \cite{Si69})

(b) The standard $d$-winding
$$S^1\to S^1\subset\R^2,\quad (\cos\varphi,\sin\varphi)\mapsto (\cos d\varphi,\sin d\varphi)$$
is approximable by embeddings if and only if $d\in\{-1,0,1\}$.

A map $S^1\to S^1\subset\R^2$ is approximable by embeddings if and only if its degree is in $\{-1,0,1\}$. (Sieklucki, 1969 \cite{Si69})

(c) An Euler path or an Euler cycle in the plane is approximable by embeddings if and only if it
 does not have transversal self-intersections. (M. Skopenkov, 2003 \cite{Sk03})

% (hence any Euler graph in the plane has an Euler cycle, approximable by embeddings).

(d) There exists an algorithm of checking whether a given simplicial map is approximable by embeddings (A. Skopenkov, 1994, M. Skopenkov, 2003
\cite{Sk94, Sk03}).
\end{proposition}

%The other proofs are more complicated.

The above examples and results show that approximability by embeddings is non-trivial to check, even for paths.
One can also show that there is no Kuratowski-type criterion.

\begin{problem}
If a polygonal line of $k$ edges is approximable by polygonal lines without self-intersections, is it then approximable by polygonal lines of $k$ edges without self-intersections?
\end{problem}

The remaining subsections of this section are not used later (except that \S\ref{s:apart} is used in Example  \ref{e:apart} and for minor remarks in \S\S\ref{s:vkcom},\ref{s:hicr}\ref{s:genr}).

\subsection{Relation to planarity of compacta}\label{s:repl}

This subsection is formally not used later, but serves as a motivation for the definition of embeddable map from \S\ref{s:exdy}.

% and so could be omitted.
% We define the decomposition of a 1-dimensional compactum into {\it an inverse limit} and show how the notion

 \begin{example}[the 2-adic van Danzig solenoid]
 Take a solid torus $T_1\subset\R^3$.
 Let $T_2\subset T_1$ be a solid torus going twice along the axis of
 the torus $T_1$.
 Analogously, take $T_3\subset T_2$ going twice along the axis of $T_2$.
 Continuing in the similar way, we obtain an infinite sequence  of solid tori $T_1\supset T_2\supset T_3\supset\dots$
 The intersection of all tori $T_i$ is
%a 1-dimensional compactum and is
called  {\it the 2-adic van Danzig solenoid}.
 \end{example}

{\bf The inverse limit} of an infinite sequence
 $K_1 \overset{\varphi_1}\leftarrow K_2 \overset{\varphi_2}\leftarrow K_3 \overset{\varphi_3}\leftarrow\dots $
 of graphs in $\R^3$ and simplicial maps between them is the compactum
 $$
 K = \{\, (x_1, x_2, \dots) \in l_2 \, : \, x_i\in K_i \text { and } \varphi_ix _ {i+1} =x _ {i} \, \}.
 $$
E.g. for the van Danzig solenoid one can take all $K_i=S^1$  and all $\varphi_i$ to be 2-windings.
Any 1-dimensional compactum can be represented as an inverse limit as above.

% Such representation  shows that any 1-dimensional compactum can be embedded into $\R^3$. It also gives
%There is a simple sufficient condition to planarity of $K$:

Clearly, {\it $K$ is planar if for each $i$ and each embedding $f_i:K_i\to \R^2$ the composition $f_i\circ \varphi_i$ is approximable by embeddings.}

\subsection{Examples that appeared in dynamical systems}\label{s:exdy}

%A partial case of the approximability problem, important

A map $f:K\to M$ is said to be {\it embeddable} in $X$ if there exists an embedding $\psi:M\to X$ for which
$\psi\circ f$ is approximable by embeddings.

Let $K$ and $M$ be wedges of $k$ and $m$ circles, respectively, and suppose that
$f$ is represented by $k$ words of $m$ letters.

 \begin{example}
(Smale) The map $S^1\vee S^1\to S^1\vee S^1$, defined by $a\mapsto aba$ and
$b\mapsto ab$ is embeddable into torus but not into plane.

(Wada--Plykin) The map $S^1\vee S^1\vee S^1\to S^1\vee S^1\vee S^1$, defined by $a\mapsto aca^{-1}$, $b\mapsto bab^{-1}$ and $c\mapsto b$
is embeddable into plane.

(Zhirov) The map $S^1\vee S^1\vee S^1\vee S^1\to S^1\vee S^1\vee S^1\vee S^1$,
defined by $a\mapsto ac$, $b\mapsto ad$, $c\mapsto bac$ and $d\mapsto c$ is
embeddable into pretzel but not into torus.
 \end{example}

\subsection{Moving graphs in the plane apart}\label{s:apart}

%{\bf Approximability by maps with disjoint images}

A map $\varphi:K\sqcup L\to \R^2$ of the disjoint union of graphs $K$ and $L$ is {\bf approximable by maps with disjoint images},
or {\bf disjoinable}, if for each $\varepsilon>0$ there is an $\varepsilon$-close to $\varphi$ map $f:K\sqcup L\to\R^2$ such that
$f(K)\cap f(L)=\emptyset$.

Clearly,

$\bullet$ if $K=L=I$ and $\varphi(K)=\varphi(L)=I\times0\subset\R^2$, then $\varphi$ is disjoinable.

$\bullet$ transversal intersection is not disjoinable.

$\bullet$ there is a non-disjoinable map $\varphi:I\sqcup I\to \R^2$ having no transversal intersections.

\begin{problem} Is a polygonal line $\varphi:I\to \R^2$ approximable by embeddings if
for each pair $I_1,I_2\subset I$ such that $I_1\cap I_2=\emptyset$ the pair of of subdivided polygonal lines $\varphi|_{I_1\sqcup I_2}$ is disjoinable?
\end{problem}

This might follow from Theorem \ref{t:minc}.a of Minc below.

Analogously to approximability by embeddings one can prove that there exists an algorithm of checking disjoinability of simplicial maps.

\begin{problem} Find a `quick' algorithm of checking disjoinability of simplicial maps, at least for $K=L=I$, i.e., for polygonal lines.
\end{problem}

\subsection{Acknowledgements}
This short survey is based on a talk at Institute of Science and Technology Austria (August, 2016),
Dept. of Computer Sciense at HSE (Moscow, October, 2016) and MIPT conference (Dolgoprudnyi, November, 2016).
I am grateful to U. Wagner for the invitation to IST, to  H. Edelsbrunner, R. Fulek, D. Repov\v s, M. Skopenkov and anonymous referee for useful discussions, and to M. Skopenkov for allowing me to use some figures.
Supported in part by the Russian Foundation for Basic Research Grant No. 15-01-06302, by Simons-IUM Fellowship
and by the D. Zimin Dynasty Foundation.

\section{Criteria for approximability by embeddings of paths}\label{s:crit}

Let us state criteria for  approximability by embeddings of a simplicial path or cycle in the plane.
These criteria assert that, in some sense, transversal self-intersections are the main obstructions to approximability by embeddings.
The criteria are stated in terms of {\it the derivative} \cite{Mi97} which we define later.

This section is not used in the rest of this paper (except for \S\ref{s:vkcom}).

\begin{theorem}[P. Minc, 1997, M. Skopenkov, 2003 \cite{Mi97, Sk03}]\label{t:minc}
(I) A polygonal line (i.e. PL map) $\varphi:I\to \R^2$ with $k$ vertices is approximable by embeddings if and only if for each $i=0,\dots,k$ the $i$-th derivative $\varphi^{(i)}$ does not have transversal self-intersections.

(S) A closed polygonal line $\varphi:S^1\to \R^2$ with $k$ vertices is approximable by embeddings if and only if for each $i=0,\dots,k$ the $i$-th derivative $\varphi^{(i)}$ neither contains transversal self-intersections nor is the standard winding of degree $d\not\in\{-1,0,1\}$.
\end{theorem}

Motivated by Theorem \ref{t:vkd} below, we conjecture that the analogues of these criteria for simplicial maps $\varphi:K\to \R^2$ of connected graphs $K\not\cong I,S^1$ are false.

\begin{figure}[h]\centering
\includegraphics[width=11.5cm]{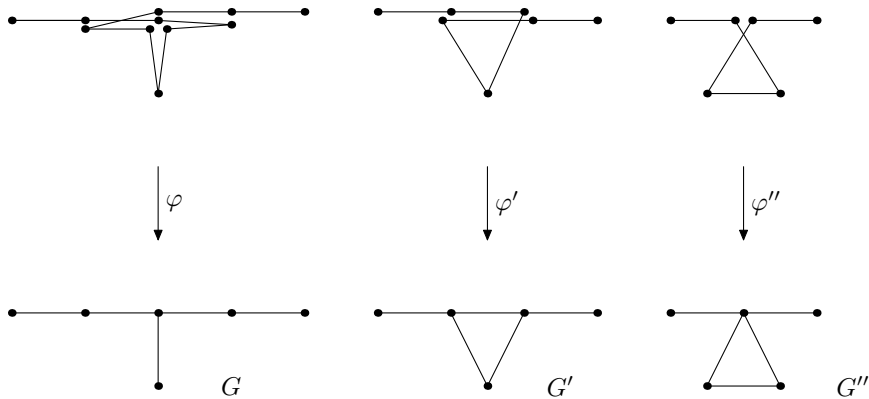}
\caption{Derivatives of graphs and paths}
\label{t-fig1}
\end{figure}

Let us define the {\bf derivative} (= the {\bf line graph}) $G'$ of a graph $G$.
 The vertex set of the graph $G'$ is the edge set of $G$.
 For an edge $a\subset G$ denote by $a'\in G'$ the corresponding vertex.
 Vertices $a'$ and $b'$ of $G'$ are joined by an edge if and only if the edges $a$ and $b$ are adjacent in $G$.

% Note that the derivatives $G'$ and $H'$ of homeomorphic but not isomorphic graphs $G$ and $H$
%are not necessarily homeomorphic.

Now let $\varphi$ be a path in the graph $G$ given by the sequence of vertices $x_1,\dots,x_k\in G$, where
 $x_i$ and $x_{i+1}$ are joined by an edge.
 Then $(x_1x_2)',\dots,(x_{k-1}x_k)'$  is a sequence of vertices of $G'$.
 In this sequence replace consecutive same vertices  by a single  vertex.
 The obtained sequence of vertices determines a path $\varphi'$ in the graph $G'$ called {\bf pre-derivative} of the path $\varphi$.

% each collection  $(x_ix_{i+1})'$, $(x_{i+1}x_{i+2})$, $\dots$, $(x_{j-1}x_{j})'$ of

A 5-od (the cone over 5 points) is a planar graph whose derivative is the Kuratowsky graph, which is not planar.
But if $G\subset\R^2$ and the path $\varphi$ does not have transversal self-intersections, then the image of the map $\varphi'$ is a subgraph $G'_\varphi\subset G'$ {\it naturally embedded in the plane} \cite{Sk03}.
The {\bf derivative} of the path $\varphi$ is defined by changing $G'$ to $G'_\varphi\subset\R^2$ and $\varphi'$ to its onto restriction $\varphi':I\to G'_\varphi$.

%(we give the construction of a natural embedding $G'_\varphi\to\R^2$ in \S2, Definition~of~$N'$).

Define the {\bf $k$-th derivative} $\varphi^{(k)}$ inductively.

For a cycle $\varphi$ the definition of the  {\bf derivative} cycle $\varphi'$ is analogous.

\begin{remark}
(a)  $\varphi'=\varphi$ for the standard $d$-winding $\varphi:S^1 \to S^1$ with $d\ne0$.

(b) $\varphi'$ is an embedding for any Euler path or cycle $\varphi$ having no transversal self-intersections.

% Thus Example~1.2 is indeed a specific case of the following theorem.

(c)  Let $\varphi:S^1\to G$ be a cycle of $k$ vertices in a graph $G$.
 Then either the domain of $\varphi^{(k)}$ is empty or $\varphi^{(k)}$ is a standard winding of a non-zero degree.
This degree can be considered as a {\bf generalization of the degree of a map $S^1\to S^1$}.
\end{remark}

\begin{corollary}[P. F. Cortese, G. Di Battista, M. Patrignani, M. Pizzonia, 2009 \cite{CBPP}, cf. \cite{CEX, AAFT}]
There are polynomial algorithms for

$\bullet$ checking approximability of paths or cycles by embeddings.

$\bullet$ finding the degree of the winding $\varphi^{(k)}$ for a cycle $\varphi:S^1\to G$ of $k$ vertices in a graph $G$.
\end{corollary}

%Such an algorithm will make criteria for approximability by embeddings easier to apply.

This result is proved in \cite{CBPP} directly, without reference to Theorem \ref{t:minc} of Minc-Skopenkov.
But apparently Theorem \ref{t:minc} was reproved, without stating it explicitly, in the course of proving this result.

See also \cite{Fu} on approximability by embeddings having a fixed isotopy class.

\begin{figure}[h]\centering
\includegraphics[width=12cm]{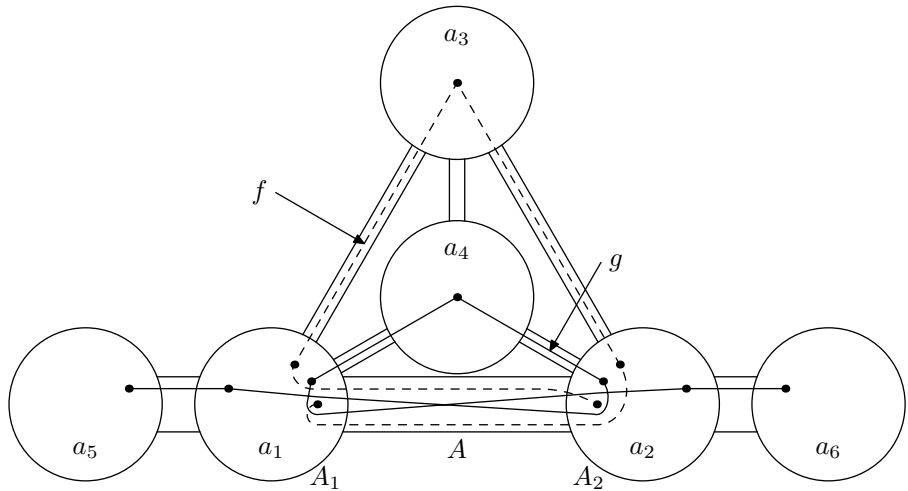}
\caption {A pair of non-disjoinable paths with disjoinable derivatives}
\label{t-fig11}
\end{figure}

\begin{example}[M. Skopenkov, 2003 \cite{Sk03}]\label{e:apart} There exists a pair of polygonal lines $\varphi,\psi:I\to \R^2$ (on the figure above a pair of paths $f,g$, close to $\varphi,\psi$, is shown),
which is not disjoinable but for which the pair $\varphi',\psi'$ of derivatives is disjoinable.
\end{example}

\section{The van Kampen obstruction}

\subsection{The Hanani-Tutte obstruction}

{\it The Hanani-Tutte obstruction for planarity of graphs} appeared in works of Hanani and Tutte in 1934, 1970.

{\bf The Hanani-Tutte obstruction} for approximability of a polygonal line $\varphi:I\to\R^2$ by embeddings is
the following necessary condition:

for each $\varepsilon>0$ there exist a subdivision $J$ of $I$ and a general position polygonal line $f:J\to\R^2$ which is $\varepsilon$-close to $\varphi$ and such that $|f(e_1)\cap f(e_2)|$ is even for each disjoint edges $e_1,e_2$ of $I$.

Note that

$\bullet$ $f(e_1)$ and $f(e_2)$ are polygonal lines, not necessarily segments.

$\bullet$ for small enough $\varepsilon>0$ we have $f(e_1)\cap f(e_2)=\emptyset$ whenever
$\varphi(e_1)\cap \varphi(e_2)=\emptyset$.

{\it The integer Hanani-Tutte obstruction} is defined analogously, the only change is that instead of `$|f(e_1)\cap f(e_2)|$ is even' we orient edges $e_1,e_2$ and require that the sum of {\it signs} of their intersection points is zero.

\subsection{Definition of the van Kampen obstruction}\label{s:vk}

The van Kampen obstruction is a reformulation of the Hanani-Tutte obstruction in an equivalent form more convenient to construct algorithms.
{\it The van Kampen obstruction for embeddability of $n$-dimensional complex in $\R^{2n}$} was invented by van Kampen around 1932.
For an elementary exposition of the case $n=1$ see e.g. \cite{Fo04}, \cite[\S1]{Sk18}, \cite[\S1]{Sk}.

Let us define the mod 2 van Kampen obstruction to approximability by embeddings of polygonal lines.
(This construction is more visual than that for of embeddability of graphs or complexes.)

 Let $\varphi:I\to \R^2$ be a simplicial path.
 Denote by $x_1,\dots,x_k$ the vertices of $I$ in the order along $I$,
 and denote the edge $x_ix_{i+1}$ by $i$.
 Let $I^*=\bigcup\limits_{i<j-1}i\times j$
 be the {\it simplicial deleted product} of given subdivision of $I$.

Paint black the edges $x_i\times j$, $j\times x_i$, and the cells $i\times j$ of $I^*$ such that $\varphi x_i\cap\varphi j=\emptyset$,  $\varphi i\cap\varphi j=\emptyset$.
Denote by $I^{*\varphi}$ the black set.

\begin{figure}[h]\centering
\includegraphics{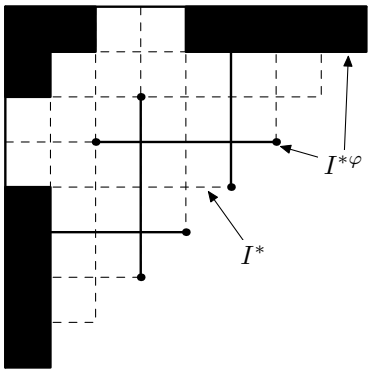}
\caption{The Van Kampen obstruction for the `T-path' (fig. \ref{putibcd}.a).}
\label{t-fig9}
\end{figure}

 Take a general position map $f:I\to\R^2$ sufficiently close to $\varphi$.
 To each cell $i\times j$ of `the table' $I^*$ put the number
 $$v_f(i\times j):=|fi\cap fj|\pmod2.$$
Cut $I^*$ along the black edges, and let $C_1,C_2,\dots,C_n$ be all the obtained components such that
 $\partial C_k\cap\partial I^*\subset I^{*\varphi}$.
 Denote
 $$v_f(C_k)=\sum\limits_{i\times j\subset C_k}v_f(i\times j).$$
 {\it The van Kampen obstruction} (with $\Z_2$-coefficients) for approximability by embeddings is the vector
 $$v(\varphi)=(\,v_f(C_1),v_f(C_2),\dots ,v_f(C_n)\,).$$
One can easily show that $v(\varphi)$ does not depend on the choice of $f$.

Thus $v(\varphi)=0$ is a necessary condition for approximability by embeddings.
It is also easy to check that $v(\varphi)\ne0$ for a polygonal line $\varphi:I\to\R^2$
 containing a transversal self-intersection.

\subsection{Completeness of the van Kampen obstruction}\label{s:vkcom}

\begin{corollary}[M. Skopenkov, 2003 \cite{Sk03}]\label{c:vki}
A polygonal line $\varphi:I\to \R^2$ is approximable by embeddings if and only if the van Kampen obstruction $v(\varphi)$ is zero.
\end{corollary}

This criterion, although more difficult to state, might be used to construct a faster algorithm than the `derivative' criterion of Theorem \ref{t:minc}.a above.

{\it The van Kampen obstruction} $v(\varphi)$ for approximability by embeddings of a simplicial map $\varphi:K\to \R^2$ of any graph $K$ is defined analogously to \S\ref{s:vk} \cite[\S1]{RS}, \cite{Sk03}.

%(in the lecture or in \cite[\S1]{RS}).

The analogue of Corollary \ref{c:vki} for {\it closed} polygonal line is false.
The standard 3-winding is a counterexample \cite[Fig. 5a]{RS}.
Using `integration' one can construct a counterexample with the image a triod \cite[Fig. 18]{FKMP}.

\begin{conjecture}\label{p:vkt} \cite[Conjecture 1.8]{RS} \cite[Conjecture 1.6]{Sk03} Is analogue of Corollary \ref{c:vki} true for the polygonal line replaced by a simplicial map $\varphi:K\to \R^2$ of any tree $K$?
\end{conjecture}

This conjecture holds when $\varphi(K)\subset\R$ \cite{Fu14, Fu16}.
A positive solution for general case is announced in \cite{FK}; some of the mathematicians who posed the conjecture do not recognize that proof as complete.

\begin{theorem}[M. Skopenkov, 2003 \cite{Sk03}]\label{t:vkd}
 Let $K$ be a graph with $k$ vertices and without vertices of degree $>3$.
 A simplicial map $\varphi:K\to S^1\subset\R^2$ is approximable by embeddings if and only if the van Kampen obstruction $v(\varphi)=0$ and $\varphi^{(k)}$ does not contain standard windings of odd degree $d\ne1$.
\end{theorem}

Here the derivative $\varphi^{(k)}$ is defined analogously to \S\ref{s:crit}.

\begin{problem} Is analogue of Theorem \ref{t:vkd} true for an arbitrary graph $K$ and $S^1$ replaced by an arbitrary graph $G$?
\end{problem}

A positive answer implies a positive solution of Conjecture \ref{p:vkt}.

One can define {\it the van Kampen obstruction for  disjoinability} (see \S\ref{s:apart})
However, Example \ref{e:apart}
shows that its triviality is not sufficient for disjoinability, even for $K=L=I$, i.e., for polygonal lines.

%the above example of M. Skopenkov
% V) {\it The van Kampen obstruction} $v(\varphi,\psi)=0$.
%The map $\Phi:\{(x,y)\in I\times I\ |\ \varphi x\ne\psi y\}\to S^1$ given by
%$\Phi(x,y)=\frac{\varphi x-\psi y}{\|\varphi x-\psi y\|}$ homotopically extends to a map $I\times I\to S^1$.

\subsection{The deleted product}\label{s:delpro}

To formulate some results and conjectures we need the important notion of the {\it deleted product}.

\begin{remark} (This remark introduces the idea of the deleted product but is formally not used later.)
Alice and Bob stand in vertices $A,B$ of the triangle $ABC$.
They can walk continuously along the edges of the triangle, so that each moment one of them is in a vertex, and the other is on the opposite edge.
Then they can exchange their positions.
\end{remark}

The {\it deleted product} $\t G$ of a graph (or a topological space) $G$ is the product of
$G$ with itself, minus the diagonal:
$$\t G:=\{(x,y)\in G\times G\ |\ x\ne y\}.$$
This is the configuration space of ordered pairs of distinct points of $G$.

\begin{figure}[h]
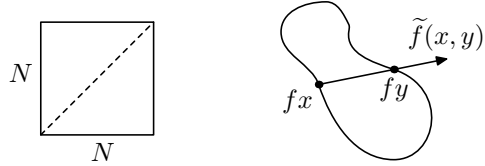
\centering
\includegraphics{4-2.eps}\qquad\qquad\includegraphics[scale=1]{4-3.eps}
\caption{The deleted product and the Gauss map}\label{}
\end{figure}

%\begin{figure}[h]\centering
%\includegraphics{4-3.eps}
%\caption{The Gauss map}\label{}
%\end{figure}

Now suppose that $f:G\to\R^m$ is an embedding (e.g. $m=2$).
Then the map $\t f:\t G\to S^{m-1}$ is well-defined by the Gauss formula
$$\t f(x,y)=\frac{f(x)-f(y)}{|f(x)-f(y)|}.$$
We have $\t f(y,x)=-\t f(x,y)$.
I.e. this map is equivariant with respect to the `exchanging factors' involution
$(x,y)\to(y,x)$ on $\t G$ and the antipodal involution on $S^{m-1}$.
Thus the existence of an equivariant map $\t G\to S^{m-1}$ is a
necessary condition for the embeddability of $G$ in $\R^m$.

\subsection{How close is a map to an embedding?}

Let $\varphi: G\to\R^2$ be a map from a (finite) graph to the plane.

Denote by $a(\varphi)$ the infimum of those $a\ge0$ for which there is a map $f: G\to\R^2$ which is $a$-close to $\varphi$ and has no self-intersections (i.e. is an embedding).

%(The distance between maps is defined as the maximum Euclidean distance between points $\varphi(x)$ and $f(x)$
%over all $x\in S^1$.)

\begin{problem}[H. Edelsbrunner]\label{p:a}  Estimate $a(\varphi)$, and construct an embedding $f$ realizing the estimation from above.
\end{problem}

%Denote $\t G:=\{(x,y)\in G\times G\ |\ x\ne y\}.$

Denote by $b(\varphi)$ the infimum of those $b\ge0$ for which there is an equivariant map $\t G\to S^1$
whose restriction to the set
$$\{ (x,y)\in G\times G\ :\ |\varphi x-\varphi y|>b\}$$
is equivariantly homotopic to the map given by the formula
$$\t\varphi(x,y):=\dfrac {\varphi(x)-\varphi(y)} {|\varphi(x)-\varphi(y)|}.$$
The number $b(\varphi)$ can be constructively estimated using {\it the van Kampen obstruction}.
Clearly, $a(\varphi)\ge b(\varphi)/2$.

\begin{problem}[appeared in a discussion with H. Edelsbrunner]\label{l:est} Is there a number $C$ such that
$a(\varphi)\le C b(\varphi)$ for each polygonal line  $\varphi:[0,1]\to\R^2$?
(So $C$ does not depend of the number of vertices of the polygonal line.)
\end{problem}

One can try to obtain an affirmative solution of the above problem
(and to construct a map realizing the estimation from above) using
proof of a criterion for approximability by embeddings \cite[Corollary 1.4]{Sk03}.

The answer for  a cycle $\varphi:S^1\to\R^2$ is presumably `no'.
However, one can possibly formulate a more elaborate conjecture using derivation of graphs
\cite[Theorems 1.3.S and 1.5]{Sk03}.

The answer for  a higher-dimensional analogue of the Problem is presumably `yes'.
One can use the proof of a criterion for approximability by embeddings \cite{RS}.
That proof  presumably gives a construction of a map realizing the estimation from above.

\section{Higher dimensional generalizations}

\subsection{Approximability by embeddings in higher dimensions}

A map $f:K\to M$ between polyhedra (=bodies of simplicial complexes) is said to be {\bf embeddable} in $\R^m$ if there exists an embedding
$\psi:M\to\R^m$ for which $\psi\circ f$ is approximable by embeddings.

\begin{theorem} (a) For each $n$ every map $f:I^n\to I^n$ is embeddable in $\R^{2n}$ (Sieklucki, 1969 \cite{Si69}).

(b) For each $n>1$ every map $f:T^n\to T^n$ between $n$-dimensional tori is embeddable in $\R^{2n}$ (Keesling-Wilson, 1985 \cite{KW}).

(c) For each $n>1$ every map $f:S^n\to S^n$ is embeddable in $\R^{2n}$ (Akhmetiev 1996, Melikhov, 2004 \cite{Ak96, Me04}).
\end{theorem}

The proof of (a) is obvious, the proofs of (b,c) are much more complicated.

%\section{General results on approximability by embeddings}

\begin{corollary}
For each fixed $m,n$ such that $2m\ge3n+3$ approximability by  embeddings of simplicial maps $K\to \R^m$ from $n$-complexes is decidable in polynomial time.
\end{corollary}

%(\cite{RS}, see

This follows by the 1998 algebraic (=combinatorial) criterion for approximability by  embeddings,
Theorem \ref{t:vkhi} below, and from the 2013 result of M. \v Cadek-M. Kr\v cal-L. Vok\v rinek on algorithmic solvability of the corresponding algebraic problem \cite{CKV}.

For each fixed $m,n$ such that $2m<3n+3$ approximability by  embeddings is NP hard by Matou\v sek-Tancer-Wagner, 2008 \cite{MTW}.
Cf. \cite[\S6]{AAFT}.

\subsection{Criteria for approximability by embeddings}\label{s:hicr}

Let $K,L$ be $n$-dimensional complexes.
{\it The van Kampen obstructions} for embeddability of $K$ into $\R^{2n}$, for approximability by embeddings of a simplicial map
$\varphi:K\to \R^{2n}$, and for approximability of a simplicial map $\varphi:K\sqcup L\to \R^{2n}$ by maps with disjoint images, are defined analogously to \S\ref{s:vk} \cite[\S4]{Sk08}, \cite[\S1]{RS}.

\begin{theorem}[E.R. Van Kampen, A. Shapiro, W.T. Wu, 1932-57]
For $n\ne2$ an $n$-complex is embeddable into $\R^{2n}$ if and only if the van Kampen obstruction for embeddability of $K$ into $\R^{2n}$ is zero.
\end{theorem}

%\cite{Sk08}]

The analogue for $n=2$ is wrong by M.~Freedman, V.~Krushkal and P.~Teichner, 1994 \cite{FKT}
(see a simpler proof in \cite[\S2.2]{AMSW}).

\begin{theorem}[D. Repov\v s and A. Skopenkov, 1998 \cite{RS}]\label{t:vkhi}
For $n>2$ a simplicial map $K\to \R^{2n}$ of an $n$-complex is approximable by embeddings if and only if the van Kampen obstruction for approximability by embeddings is zero.
\end{theorem}

The analogue for $n=1,2$ is wrong.

%Analogous (but simpler) result for disjoinability was essentially proved
%by S.~Spie\. z and H.~Toru\'nczyk in 1991.

%\section{Completeness of On the deleted product obstruction}

The deleted product $\t K$ is defined in \S\ref{s:delpro} above.

%\cite{Sk08}]

\begin{theorem}[C. Weber, 1967]
Assume that $2m\ge3n+3$ and $K$ is an $n$-complex.
There is an embedding $K\to \R^m$ of  if and only if there is an equivariant map $\t K\to S^{m-1}$.
\end{theorem}

\begin{theorem}[D. Repov\v s and A. Skopenkov, 1998 \cite{RS}]
Assume that $2m\ge3n+3$ and $\varphi:K\to \R^m$ is a simplicial map of an $n$-complex $K$.
The map $\varphi$ is approximable by  embeddings if and only if there is
an equivariant map $\t K\to S^{m-1}$ whose restriction to the set $\{ (x,y)\in K\times K: \varphi x\ne\varphi y\}$ is equivariantly homotopic to the map given by the formula $\t\varphi(x,y)=\dfrac {\varphi(x)-\varphi(y)} {|\varphi(x)-\varphi(y)|}$.
 \end{theorem}

Analogous simpler result for disjoinability was proved by S.~Spie\. z and H.~Toru\'nczyk in 1991 \cite{ST}.
They were motivated by the problem of approximability of {\it any map} from {\it given compacta} by maps having disjoint images.
A different approach to this problem was proposed at the same time by A. N. Dranishnikov, D. Repov\v s and E. V. Shchepin, see \cite{Dr00} and references therein.

%%The different approach yielded stronger results on this problem.
% (comparatively to \cite{ST91}).
%%But it is not clear how to apply that approach to approximability of a
%%{\it given map} by maps having disjoint images.

\subsection{Generalizations to $r$-fold points}\label{s:genr}

\begin{problem}[Gromov, 2010 \cite{Gr10}]
Is it correct that if $r$ is not a prime power, then for each compact subset $K$ of $\R^m$ for some $m$, having Lebesgue dimension $\dim K=(r-1)n$,
 there is a continuous map $K\to\R^{nr}$ each of whose point preimages contains less than $r$ points?
\end{problem}

%The answer is `yes'??? for polyhedra $K$ and each $n\ge2$ \cite{MW, AMSW}.

We conjecture that the answer is `no'.

Embeddability of compacta is related to approximability of maps from polyhedra by embeddings.
Analogously, Gromov's problem is related to approximability of maps from polyhedra by maps without $r$-fold points.

Let $K$ be a finite simplicial complex.
A map  $f\colon K\to \R^m$ is an {\bf almost $r$-embedding}
if $f(\sigma_1)\cap \dots \cap f(\sigma_r)=\emptyset$ whenever $\sigma_1,\dots,\sigma_r$ are pairwise disjoint simplices of $K$.

\begin{theorem}[I. Mabillard and U. Wagner, 2015, S. Avvakumov, I. Mabillard, A. Skopenkov and U. Wagner, 2015 \cite{MW, AMSW}]
If $n\ge2$ and $r$ is not a prime power, then any $n(r-1)$-complex is almost $r$-embeddable in $\R^{nr}$.
\end{theorem}

{\it The van Kampen obstruction} for almost $r$-embeddability of an $n(r-1)$-complex to $\R^{nr}$
is defined analogously to \S\ref{s:vk} above \cite{MW}.

\begin{theorem}[I. Mabillard and U. Wagner, 2015, S. Avvakumov, I. Mabillard, A. Skopenkov and U. Wagner, 2015 \cite{MW, AMSW}]
Suppose that $n,r\ge2$, $n+r\ge5$.
An $n(r-1)$-complex is almost $r$-embeddable in $\R^{nr}$ if and only if the van Kampen obstruction is zero.
\end{theorem}

\begin{problem} Are analogous results true for approximability by almost $r$-embeddings, or by maps $f:K_1\sqcup\ldots\sqcup K_r\to\R^{nr}$ such that $f(K_1)\cap\ldots\cap f(K_r)=\emptyset$?
\end{problem}

%for $r$-disjoinability, or

%\small

\end{document}